\documentclass[a4paper,12pt]{report}
\usepackage{inputenc}
\usepackage[english]{babel}
\usepackage[tbtags]{amsmath}
\usepackage{amsfonts,amssymb,mathrsfs,amscd}
\usepackage{graphics}
\usepackage{wrapfig}
\usepackage[dvips]{graphicx}

\date{}
\newtheorem{df}{Definition}[section]
\newtheorem{Th}{Theorem}[section]
\newtheorem{Lm}[Th]{Lemma}

\begin{document}
\renewcommand{\baselinestretch}{0.96}
\renewcommand{\thesection}{\arabic{section}}
\renewcommand{\theequation}{\thesection.\arabic{equation}}
\csname @addtoreset\endcsname{equation}{section}
 \large
\begin{center}
\textbf {S.\,F.~Lukomskii \\
Step refinable functions and orthogonal MRA on $p$-adic Vilenkin
groups
 \footnote[1]{This research was carried out with the financial
support the Russian Foundation for Basic Research (grant
no.~10-01-00097).}\\
(Russia, Saratov)\\
lukomskiisf@info.sgu.ru\\}
\end{center}

\vskip1cm

\qquad \parbox{13cm}{{\bf Abstract.}
  We find the necessary and
 sufficient conditions  for refinable step function under which this
 function  generates an orthogonal MRA in the $L_2(\mathfrak G)$
 -spaces on Vilenkin groups $\mathfrak G$. We consider a class of refinable step
functions for which the mask $m_0(\chi)$ is constant on cosets
$\mathfrak G_{-1}^\bot$ and its modulus $|m_0(\chi)|$ takes two
values only: 0 and 1. We will prove that any refinable step
function $\varphi$ from this class that generates an orthogonal
MRA on $p$-adic Vilenkin group $\mathfrak G$ has Fourier transform
with condition ${\rm supp}\,\hat \varphi(\chi)\subset \mathfrak
G_{p-2}^\bot$. We show the sharpness of this result too.}

 \vskip1cm


\section*{Introduction}
 Foundations for wavelet analysis theory on locally compact groups
 have been lay in the monograph [1]. In articles [2-4] first
 examples of orthogonal wavelets on the dyadic Cantor group are
 constructed and their properties are studied. The general scheme
 for the construction of wavelets  is based on the notion of
 multiresolution analysis (MRA in the sequel) introduced by Y.
 Meyer and S. Mallat [5, 6]. Yu.Farkov [7-12] found necessary and
 sufficient conditions  for a refinable function under which this
 function  generates an orthogonal MRA in the $L_2(\mathfrak G)$
 -spaces on the Vilenkin group $\mathfrak G$. These conditions use the
 Strang-Fix and the modified Cohen  properties. In [10]  this
 construction  to the p = 3 case in a concrete fashion  are given.
 In [11], some algorithms for constructing orthogonal and biorthogonal compactly
 supported wavelets on Vilenkin groups are suggested. In [7-11] two
 types of orthogonal wavelets examples are constructed: step
 functions and sums of Vilenkin  series.

In these examples all step refinable functions have a support
${\rm supp}\,\hat \varphi(\chi)\subset \mathfrak
G_1^\bot=\mathfrak G_0^\bot{\cal A}$ where $\mathfrak G^{\perp}_0$
is  the unit ball in the character group and ${\cal A}$ is a
dilation operator. Therefore there is an assumption that a step
refinable function which generates an orthogonal MRA on Vilenkin
group $\mathfrak G$ has a Fourier transform with support ${\rm
supp}\,\hat \varphi(\chi)\subset \mathfrak G_1^\bot$ . We will
prove that it is not true. We consider a class of refinable step
functions for which the mask $m_0(\chi)$ is constant on cosets
$\mathfrak G_{-1}^\bot$ and its modulus $|m_0(\chi)|$ takes two
values only: 0 and 1. We will prove that any refinable step
function $\varphi$ from this class that generates an orthogonal
MRA on $p$-adic Vilenkin group $\mathfrak G$ has Fourier transform
with condition ${\rm supp}\,\hat \varphi(\chi)\subset \mathfrak
G_{p-2}^\bot$. We show the sharpness of this result too.

We should note that in the p-adic analysis, the situation is
different. S. Albeverio, S. Evdokimov, M. Skopina [12] proved,
that if a refinable step function $\varphi$ generates an
orthogonal $p$-adic MRA, then $\hat \varphi(\chi)\subset \mathfrak
G_0^\bot$.

\section{Preliminaries}
 We will consider the Velenkin group as a  locally compact zero-dimensional Abelian group
 with additional condition $p_ng_n=0$. Therefore we start with some basic notions
  and facts related to analysis on zero-dimensional groups. A topological group in which the connected component of~0 is~0 is
usually referred to as a~\textit{zero-dimensional group}. If
a~separable locally compact group $(G,\dot + )$ is
zero-dimensional, then the topology on it can be generated by
means of a~descending sequence of subgroups.

The converse statement holds for all topological groups
(see~\cite[Ch.~1, \S\,3]{13}). So, for a~locally compact group, we
are going to say `zero-dimensional group' instead of saying
`a~group with topology generated by a~sequence subgroups'.

Let $(G,\dot + )$~be a~locally compact zero-dimensional Abelian
group with the topology generated by a~countable system of open
subgroups
$$
\cdots\supset G_{-n}\supset\cdots\supset G_{-1}\supset G_0\supset
G_1\supset\cdots\supset G_n\supset\cdots
$$
where
$$
\bigcup_{n=-\infty}^{+\infty}G_n= G,\quad
 \quad \bigcap_{n=-\infty}^{+\infty}G_n=\{0\}
$$
(0~is the null element in the group~$ G$). Given any fixed
$N\in\mathbb Z$, the subgroup~$G_N$ is a~compact Abelian group
with respect to the same operation~$\dot+ $ under the topology
generated by the system of subgroups
$$
 G_N\supset G_{N+1}\supset\cdots\supset
G_n\supset\cdots.
$$
As each subgroup $G_n$ is compact, it follows that each quotient
group~$G_n/G_{n+1}$ is finite (say, of order~$p_n$). We may always
assume that all~$p_n$ are prime numbers. We will name such chain
as \it basic chain. \rm  In this case, a~base of the topology is
formed by all possible cosets~$G_n\dot + g$, $g\in G$.

We further define the numbers $(\mathfrak
m_n)_{n=-\infty}^{+\infty}$ as follows:
$$
\mathfrak m_0=1,\qquad \mathfrak m_{n+1}=\mathfrak m_n\cdot p_n.
$$
Clearly, for $n\ge1$,
$$
\mathfrak m_n=p_0p_1\cdots p_{n-1}, \qquad \mathfrak
m_{-n}=\frac1{p_{-1}p_{-2}\cdots p_{-n}}.
$$

The collection of all such cosets $G_n\dot+ g$, $n\in\mathbb Z$,
along with the empty set form the semiring~${\mathscr K}$. On each
coset $ G_n\dot+ g$ we define the measure~$\mu$ by $\mu( G_n\dot+
g)=\mu G_n=1/{m_n}$. So,~if $n\in\mathbb Z$ and~$p_n=p$, we have
$\mu G_n\cdot\mu G_{-n}=1$. The measure~$\mu$ can be extended from
the semiring~${\mathscr K}$ onto the~$\sigma$-algebra (for
example, by using Carath\-eodory's extension). This gives the
translation invariant measure~$\mu$, which agrees on the Borel
sets with the Haar measure on~$G$. Further, let
$\smash[b]{\displaystyle\int_{G}f(x)\,d\mu(x)}$~be the absolutely
convergent integral of the measure~$\mu$.

Given an $n\in\mathbb Z$, take an element $g_n\in G_n\setminus
G_{n+1}$ and fix~it. Then any $x\in G$ has a~unique representation
of the form
\begin{equation}
\label{eq1.1} x=\sum_{n=-\infty}^{+\infty}a_ng_n, \qquad
a_n=\overline{0,p_n-1}.
\end{equation}
The sum~\eqref{eq1.1} contain finite number of terms with negative
subscripts, that~is,
\begin{equation}
\label{eq1.2} x=\sum_{n=m}^{+\infty}a_ng_n, \qquad
a_n=\overline{0,p_n-1}, \quad a_m\ne 0.
\end{equation}
We will name system $(g_n)_{n\in \mathbb Z}$ as {\it a basic
system}.

Classical examples of zero-dimensional groups are Vilenkin groups
and groups of $p$-adic numbers (see~\cite[Ch.~1, \S\,2]{13}).
 A direct sum of cyclic groups $Z(p_k)$ of order~$p_k$,
 $k\in\mathbb Z$, is called a~\textit{Vilenkin group}. This means
 that the elements of a~Vilenkin group are infinite sequences
 $x=(x_k)_{k=-\infty}^{+\infty}$ such that:
 \begin{itemize}
 \item[1)] $x_k=\overline{0,p_k-1}$; \item[2)] only a~finite number
 of~$x_k$ with negative subscripts are different from zero;
 \item[3)] the group operation~$\dot+ $ is the coordinate-wise
 addition modulo~$p_k$, that~is,
 $$
 x\dot+ y=(x_k\dot+ y_k), \qquad x_k\dot+ y_k=(x_k+y_k)\ \
 \operatorname{mod}p_k.
 $$
 \end{itemize}
 A topology on such group is generated by the chain of subgroups
 $$
  G_n=\bigl\{x\in G:x=(\dots,0,0,\dots,0,x_n,x_{n+1},\dots),\
   x_\nu=\overline{0,p_\nu-1},\ \nu\ge n\bigr\}.
 $$
 The elements $g_n=(\dots,0,0,1,0,0,\dots)$ form a basic system.
 From definition of the operation $\dot+$ we have $p_ng_n=0$.
 Therefore we will name a zero-dimensional group $(G,\dot+)$ with the condition $p_ng_n=0$ as Vilenkin group.

The group~$\mathbb Q_p$ of all $p$-adic numbers ($p$~is a prime
number) also consists of sequences
$x=(x_k)_{k=-\infty}^{+\infty}$, $x_k=\overline{0,p-1}$, only
a~finite number of~$x_k$ with negative subscripts being different
from zero. However, the group operation in~$\mathbb Q_p$ is
defined differently. Namely, given elements
$$
x=(\dots,0,\dots,0,x_N,x_{N+1},\dots) \; \text{and} \;
y=(\dots,0,\dots,0,y_N,y_{N+1},\dots)\in\mathbb Q_p,
$$
we again add them coordinate-wise, but whereas in a~Vilenkin group
$x_n\dot+ y_n=(x_n+y_n)\ \operatorname{mod}p$  (that~is, a~1 is
not carried to the next $(n+1)$th position), the corresponding
$p$-adic summation has the property that the 1~occuring as
a~result of the addition of $x_n+y_n$ is carried to the next
$(n+1)$th position. We endow the group~$ \mathbb Q_p$ with the
topology generated by the same system of subgroups~$G_n$ as for
a~Vilenkin group. Similarly, as~a~$(g_n)$, we may again take the
same sequence.

By $X$ denote  the collection
 of the characters of a~group $(G,\dot+ )$; it is
a~group with respect to multiplication too. Also let
$G_n^\bot=\{\chi\in X:\forall\,x\in G_n\  , \chi(x)=1\}$ be the
annihilator of the group~$G_n$. Each annihilator~$ G_n^\bot$ is
a~group with respect to multiplication, and the subgroups~$
G_n^\bot$ form an~increa\-sing sequence
\begin{equation}
\label{eq1.3} \cdots\subset G_{-n}^\bot\subset\cdots\subset
G_0^\bot \subset G_1^\bot\subset\cdots\subset
G_n^\bot\subset\cdots
\end{equation}
with
$$
\bigcup_{n=-\infty}^{+\infty} G_n^\bot=X \quad\text{and} \quad
\bigcap_{n=-\infty}^{+\infty} G_n^\bot=\{1\},
$$
the quotient group $ G_{n+1}^\bot/ G_n^\bot$ having order~$p_n$.
The group of characters~$X$ may be equipped with the topology
using the chain of subgroups~\eqref{eq1.3}, the family of the
cosets $ G_n^\bot\cdot\chi$, $\chi\in X$, being taken as~a~base of
the topology. The collection of such cosets, along with the empty
set, forms the~semiring~${\mathscr X}$. Given a~coset $
G_n^\bot\cdot\chi$, we define a~measure~$\nu$ on it by $\nu(
G_n^\bot\cdot\chi)=\nu( G_n^\bot)= m_n$ (so that always $\mu(
G_n)\nu( G_n^\bot)=1$). The measure~$\nu$ can be extended onto the
$\sigma$-algebra of measurable sets in the standard way. One then
forms the absolutely convergent integral
$\displaystyle\int_XF(\chi)\,d\nu(\chi)$ of this measure.

The value~$\chi(g)$ of the character~$\chi$ at an element $g\in G$
will be denoted by~$(\chi,g)$. The Fourier transform~$\widehat f$
of an~$f\in L_2( G)$  is defined~as follows
$$
\widehat f(\chi)=\int_{ G}f(x)\overline{(\chi,x)}\,d\mu(x)=
\lim_{n\to+\infty}\int_{ G_{-n}}f(x)\overline{(\chi,x)}\,d\mu(x),
$$
the limit being in the norm of $L_2(X)$. For any~$f\in L_2(G)$,
the inversion formula is valid
$$
f(x)=\int_X\widehat f(\chi)(\chi,x)\,d\nu(\chi)
=\lim_{n\to+\infty}\int_{ G_n^\bot}\widehat
f(\chi)(\chi,x)\,d\nu(\chi);
$$
here the limit also signifies the convergence in the norm of~$L_2(
G)$. If $f,g\in L_2( G)$ then the Plancherel formula is valid
$$
\int_{ G}f(x)\overline{g(x)}\,d\mu(x)= \int_X\widehat
f(\chi)\overline{\widehat g(\chi)}\,d\nu(\chi).
$$
\goodbreak

Endowed with this topology, the group of characters~$X$ is
a~zero-dimensional locally compact group; there is, however,
a~dual situation: every element $x\in G$ is a~character of the
group~$X$, and~$ G_n$ is the annihilator of the group~$ G_n^\bot$.

The union of disjoint sets $E_j$ we will denote by $\bigsqcup
E_j$.

 \section{Rademacher functions and dilation operator}
 In this section we will consider zero-dimensional groups with
 condition $p_n=p$ for any $n\in \mathbb Z$.
  In this case we define the mapping ${\cal
 A}\colon G\to G$ by
 ${\cal A}x:=\sum_{n=-\infty}^{+\infty}a_ng_{n-1}$, where
 $x=\sum_{n=-\infty}^{+\infty}a_ng_n\in G$.  The mapping~${\cal A}$ is called
 a  dilation operator if~${\cal A}(x\dot+ y)={\cal A}x\dot + {\cal A}y$ for all
 $x,y\in G$. By definition, put $(\chi {\cal A},x)=(\chi, {\cal
 A}x)$.
  A character $r_n\in  G_{n+1}^{\bot}\backslash G_n^{\bot}$
  is called the Rademacher function. Let us denote
  $$
    H_0=\{h\in G: h=a_{-1}g_{-1}\dot+a_{-2}g_{-2}\dot+\dots \dot+ a_{-s}g_{-s}, s\in \mathbb
    N\},
  $$
  $$
    H_0^{(s)}=\{h\in G: h=a_{-1}g_{-1}\dot+a_{-2}g_{-2}\dot+\dots \dot+ a_{-s}g_{-s}
    \},s\in \mathbb N.
  $$
  The set $H_0$ is an analog of the set $\mathbb N$.
  \begin{Lm}
 For any zero-dimensional group\\
 1) $\int\limits_{G_0^\bot}(\chi,x)\,d\nu(\chi)={\bf 1}_{G_0}(x)$,
 2) $\int\limits_{G_0}(\chi,x)\,d\mu(x)={\bf 1}_{G_0^\bot}(\chi)$.\\
 \end{Lm}
 The first equation it was proved in [14], the second equation  is dual to first.
 \begin{Lm}
  If $p_n=p$ for any $n\in \mathbb Z$ and the mapping ${\cal A}$ is additive then \\
  1) $\int\limits_{G_n^\bot}(\chi,x)\,d\nu(\chi)=p^n{\bf
  1}_{G_n}(x)$,\\
  2) $\int\limits_{G_n}(\chi,x)\,d\mu(x)=\frac{1}{p^n}{\bf
  1}_{G_n^\bot}(\chi)$.
  \end{Lm}
 {\bf Proof.} First we prove the equation 1). Using equations
 $$
  \int\limits_Xf(\chi{\cal
  A})\,d\nu(\chi)=p\int\limits_Xf(\chi)\,d\nu(\chi),\;\;\;{\bf
  1}_{G_n^\bot}(x)={\bf 1}_{G_0}({\cal A}^nx),
 $$
 and Lemma 2.1 we have
$$
\int\limits_{G_n^\bot}(\chi,x)\,d\nu(\chi)=\int\limits_X{\bf
1}_{G_n^\bot}(\chi)(\chi,x)\,d\nu(\chi)=p^n\int\limits_X(\chi{\cal
A}^n,x){\bf 1}_{G_n^\bot}(\chi{\cal A}^n)\,d\nu(\chi)=
$$
$$
=p^n\int\limits_X(\chi,{\cal A}^nx){\bf
1}_{G_0^\bot}(\chi)\,d\nu(\chi)=p^n{\bf 1}_{G_0}({\cal
A}^nx)=p^n{\bf 1}_{G_n}(x).
$$
The second equation is proved by analogy.  $\square$

\begin{Lm}
Let  $\chi_{n,s}=r_n^{\alpha_n}r_{n+1}^{\alpha_{n+1}}\dots
r_{n+s}^{\alpha_{n+s}}$ be a character does not belong to
$G_n^\bot$. Then
$$
\int\limits_{G_n^\bot\chi_{n,s}}(\chi,x)\,d\nu(\chi)=p^n(\chi_{n,s},x){\bf
1}_{G_n}(x).
$$
\end{Lm}
{\bf Proof.} By analogy with previously we have
$$
\int\limits_{G_n^\bot\chi_{n,s}}(\chi,x)\,d\nu(\chi)=\int\limits_X{\bf
1}_{G_n^\bot}(\chi)(\chi_{n,s}\chi,x)\,d\nu(\chi)=
$$
$$
\int\limits_{G_n^\bot}(\chi_{n,s},x)(\chi,x)\,d\nu(\chi)=p^n(\chi_{n,s},x)
{\bf 1}_{G_n}(x). \;\;\square
$$
\begin{Lm}
Let
$h_{n,s}=a_{n-1}g_{n-1}\dot+a_{n-2}g_{n-2}\dot+\dots\dot+a_{n-s}g_{n-s}\notin
G_n$. Then
$$
\int\limits_{G_n\dot+h_{n,s}}(\chi,x)\,d\mu(x)=\frac{1}{p^n}(\chi,h_{n,s}){\bf
1}_{G_n^\bot}(\chi).
$$
\end{Lm}
 This lemma is dual to lemma 2.3.
 \begin{df}
Let $M,N\in\mathbb N$.
 Denote by  ${\mathfrak D}_M(G_{-N})$ the set of step-functions
 $f\in L_2(G)$ such that 1)${\rm supp}\,f\subset G_{-N}$, and 2)
 $f$ is constant on cosets $G_M$. Similarly is defined ${\mathfrak
 D}_{-N}(G_{M}^\bot)$.
 \end{df}
\begin{Lm}
 Let $M,N\in\mathbb N$. $f\in \mathfrak D_M(G_{-N})$ if and only if $\hat f\in \mathfrak
 D_{-N}(G_M^\bot)$.
\end{Lm}
{\bf Proof.} 1) Let  $f$ be a constant on cosets $G_M\dot+g$ and
${\rm supp}\,f\subset G_{-N}$. Let us show that ${\rm supp}\,\hat
f\subset G_M^\bot$. Let $\chi\notin G_M^\bot$. Then
 $$
 \hat
 f(\chi)=\int\limits_Gf(x)\overline{(\chi,x)}\,d\mu(x)=\int\limits_{G_{-N}}f(x)\overline{(\chi,x)}\,d\mu(x)=
 $$
 $$
 =\sum_{h_{M,N}\in H_M^N}\int\limits_{G_M\dot+h_{M,N}}f(x)\overline{(\chi,x)}\,d\mu(x),
 $$
  where
 $$
   H_M^N=\{h_{M,N}=a_{M-1}g_{M-1}\dot+a_{M-2}g_{M-2}\dot+\dots\dot+a_{-N}g_{-N}\}.
 $$
  By lemma 2.4
$$
\hat f(\chi)=\sum
f(G_M\dot+h_{M,N})\int\limits_{G_M\dot+h_{M,N}}\overline{(\chi,x)}\,d\mu(x)=
$$
$$
=\sum f(G_M\dot+h_{M,N})\frac{1}{p^M}\overline{(\chi,h_{M,N})}{\bf
1}_{G_M^\bot}(\chi)=0.
$$
Now we will show that $\hat f$ is constant on cosets
$G_{-N}^\bot\zeta$. Indeed let $\chi\in G_{-N}^\bot\zeta$ and
$\zeta=r_{-N}^{\alpha_{-N}}r_{-N+1}^{\alpha_{-N+1}}\dots
r_{-N+s}^{\alpha_{-N+s}}$. Then $\chi=\chi_{-N}\zeta$ where
$\chi_{-N}\in G_{-N}^\bot$. Therefore
$$
\hat
f(\chi)=\int\limits_{G_{-N}}f(x)\overline{(\chi,x)}\,d\mu(x)=\int\limits_{G_{-N}}f(x)\overline{(\chi_{-N}\zeta,x)}\,d\mu(x)=
\int\limits_{G_{-N}}f(x)\overline{(\zeta,x)}\,d\mu(x).
$$
It means that $\hat f(\chi)$ depends only on $\zeta$. The first
part is proved.  The second part is proved similarly. $\square$
\begin{Lm}
Let $\varphi\in L_2(G)$. The system $(\varphi(x\dot-h))_{h\in
H_0}$ is orthonormal if and only if the system
$\left(p^{\frac{n}{2}}\varphi({\cal A}^nx\dot-h)\right)_{h\in
H_0}$ is orthonormal.
\end{Lm}
{\bf Proof}. This lemma follows from the equation
$$
\int\limits_Gp^{\frac{n}{2}}\varphi({\cal
A}^nx\dot-h)p^{\frac{n}{2}}\overline{\varphi({\cal
A}^nx\dot-g)}\,d\mu=\int\limits_G\varphi(x\dot-h)\overline{\varphi(x\dot-g)}\,d\mu.\;\;\square
$$
\section{MRA on Vilenkin groups}
 In what follows we will consider groups $G$ for which $p_n=p$
  and $pg_n=0$ for any $n\in \mathbb Z$. We now that it is a
  Vilenkin group. We will denote a Vilenkin group as $\mathfrak G$.
  In this group  we can chouse Rademacher functions
  in various ways.
  We define Rademacher functions by the equation
  $$
  \left(r_n,\sum_{k\in\mathbb Z}a_kg_k\right)=\exp\left(\frac{2\pi i}{p}a_n\right).
  $$
  In this case
  $$
  (r_n,g_k)=\exp\left(\frac{2\pi i}{p}\delta_{nk}\right).
  $$
 Our main objective is to
 find a refinable step-function that generates an orthogonal MRA on Vilenkin group.
 \begin{df}
  A family of closed subspaces $V_n$, $n\in\mathbb Z$,
 is said to be a~multi\-resolution analysis  of~$L_2(\mathfrak G)$
 if the following axioms are satisfied:
 \begin{itemize}
 \item[A1)] $V_n\subset V_{n+1}$;
 \item[A2)] ${\vrule width0pt
 depth0pt height11pt} \overline{\bigcup_{n\in\mathbb
 Z}V_n}=L_2(\mathfrak G)$ and $\bigcap_{n\in\mathbb Z}V_n=\{0\}$;
 \item[A3)] $f(x)\in V_n$  $\Longleftrightarrow$ \ $f({\cal A} x)\in V_{n+1}$ (${\cal A}$~is a~dilation
 operator);
 \item[A4)] $f(x)\in V_0$ \ $\Longrightarrow$ \
 $f(x\dot - h)\in V_0$ for all $h\in H_0$; ($H_0$ is analog of $\mathbb
 Z$).
  \item[A5)] there exists
 a~function $\varphi\in L_2(\mathfrak G)$ such that the system
 $(\varphi(x\dot - h))_{h\in H_0}$ is an orthonormal basis
 for~$V_0$.
\end{itemize}

 A function~$\varphi$ occurring in axiom~A5 is called
a~\textit{scaling function}.
\end{df}

%Using an MRA, we shall build functions $\psi_\nu$,
%$\nu=\overline{1,p-1}$, whose dilations and translations
%$\psi_\nu(A^jx\dminus h)$ form an~orthogonal basis
%for~$L_2(\mathfrak G)$.

 Next we will follow the conventional approach. Let
 $\varphi(x)\,{\in}\, L_2(\mathfrak G)$, and suppose that
 $(\varphi(x\dot -\nobreak h))_{h\in H_0}$ is an~orthonormal
 system in~$L_2(\mathfrak G)$. With the function~$\varphi$ and the
 dilation operator~${\cal A}$, we define the linear subspaces
 $L_j=(\varphi({\cal A}^jx\dot - h))_{h\in H_0}$ and
 closed subspaces $V_j=\overline{L_j}$. It is evident that the functions $p^{\frac{n}{2}}\varphi({\cal A}x \dot-h)_{h\in H_0}$ form
 an orthonormal basis for $V_n$, $n\in \mathbb Z$. Therefore the axiom A4 is fulfilled.  If subspaces $V_j$ form
 a~MRA, then the function~$\varphi$ is said to \textit{generate}
 an~MRA in~$L_2(\mathfrak G)$. If a function $\varphi$ generates an MRA, then we obtain from the axiom A1
\begin{equation}
  \label{eq3.1}
  \varphi(x)=\sum_{h\in H_0}\beta_h\varphi({\cal
  A}x\dot-h)\;\;\left(\sum|\beta_h|^2<+\infty\right).
 \end{equation}
  Therefore we will look up a~function
 $\varphi\in L_2(\mathfrak G)$, which generates an~MRA
 in~$L_2(\mathfrak G)$, as a~solution of the refinement
 equation (\ref{eq3.1}), A solution of refinement equation (\ref{eq3.1}) is called a {\it refinable function}.
 \begin{Lm}
Let $\varphi \in \mathfrak D_M(\mathfrak G_{-N})$ be a solution of
(\ref{eq3.1}). Then
\begin{equation} \label{eq3.2}
\varphi(x)=\sum_{h\in H_0^{(N+1)}}\beta_h\varphi({\cal A}x\dot-h)
 \end{equation}
\end{Lm}

{\bf Proof.} Let us write $\varphi(x)$ in the form
\begin{equation}                                      \label{eq3.3}
\varphi(x)=\sum_{h\in H_0^{(N+1)}}\beta_h\varphi({\cal
A}x\dot-h)+\sum_{h\notin H_0^{(N+1)}}\beta_h\varphi({\cal
A}x\dot-h).
 \end{equation}
If $x\in \mathfrak G_{-N}$, then ${\cal A}x\in \mathfrak
G_{-N-1}$. Therefore ${\cal
A}x=b_{-N-1}g_{-N-1}\dot+b_{-N}g_{-N}\dot+\dots$. If  $h\notin
H_0^{(N+1)}$, then
$$
h=a_{-1}g_{-1}\dot+\dots\dot+a_{-N-1}g_{-N-1}\dot+a_{-N-2}g_{-N-2}\dot+\dots\dot+a_{-N-s}g_{-N-s},
$$
and  $a_{-N-2}g_{-N-2}\dot+\dots\dot+a_{-N-s}g_{-N-s}\ne 0$.
 Hence  ${\cal A}x\dot-h\notin H_0^{(N+1)}$ and
$\varphi({\cal A}x\dot-h)=0$. This means that
$$
\sum_{h\notin H_0^{(N+1)}}\beta_h\varphi({\cal A}x\dot-h)=0
$$
when $x\in \mathfrak G_{-N}$.

Let  $x\notin \mathfrak G_{-N}$. Then  $\varphi(x)=0$ and ${\cal
A}x\notin \mathfrak G_{-N-1}$. Hence
$$
{\cal
A}x=\sum_{k=-N-s}^{-N-2}b_kg_k\dot+\sum_{k=-N-1}^{+\infty}b_kg_k.
$$
If  $h\in H_0^{(N+1)}$, then
$h=a_{-1}g_{-1}\dot+\dots\dot+a_{-N}g_{-N}\dot+a_{-N-1}g_{-N-1}$,
 and consequently  ${\cal A}x\dot-h\notin \mathfrak G_{-N-1}$.
 Therefore
$$
\sum_{h\in H_0^{(N+1)}}\beta_h\varphi({\cal A}x\dot-h)=0.
$$
 Using equation (\ref{eq3.3}) we obtain finally
 $$
  \sum_{h\notin H_0^{(N+1)}}\beta_h\varphi({\cal A}x\dot-h)=0,
 $$
 and lemma is proved. $\square$

 \begin{Th}
 Let $\varphi\in \mathfrak D_M(\mathfrak G_{-N})$ and let
 $(\varphi(x\dot-h))_{h\in H_0}$ be an orthonormal system.
 $V_n\subset V_{n+1}$ if and only if  the function $\varphi(x)$ is
 a solution of refinement equation (\ref{eq3.2}).
  \end{Th}
 {\bf Proof.}  First we prove that $V_n\subset V_{n+1}$ if and only if
 $V_0\subset V_1$. Indeed, let $V_0\subset V_1$ and   $f\in V_n$. Then
$$
f(x)=\sum_hc_h\varphi({\cal A}^nx\dot-h)\Rightarrow f({\cal
A}^{-n}x)=\sum_hc_h\varphi(x\dot-h)\Rightarrow f({\cal
A}^{-n}x)\in V_0\Rightarrow $$
$$
 \Rightarrow
f({\cal A}^{-n}x)\in V_1\Rightarrow  f({\cal
A}^{-n}x)=\sum_h\gamma_h\varphi({\cal A}x\dot-h)\Rightarrow
$$
$$
 \Rightarrow f(x)=\sum_h\gamma_h\varphi({\cal A}^{n+1}x\dot-h)\Rightarrow
f\in V_{n+1}.
$$
So we have, $V_n\subset V_{n+1}$. The converse is proved by
analogy.

 Now we prove that $V_0\subset V_{1}$ if and only if  the function $\varphi(x)$ is
 a solution of the  refinement equation (\ref{eq3.2}).
The necessity is evident. Let  $\varphi$ be a
 solution of (\ref{eq3.2}). We take $f\in{\rm
 span}(\varphi(x\dot-h))_{h\in H_0}$. Then
 $$
 f(x)=\sum_{\tilde h\in H_0^{(m)}}c_{\tilde h}\varphi(x\dot-\tilde
 h)
 $$
 for some  $m\in\mathbb N$.

  Since  $\varphi$ is a solution of (\ref{eq3.2}) then we can write $f$ in the form
$$
f(x)=\sum_{\tilde h\in H_0^{(m)}}c_{\tilde h}\sum_{h\in
H_0^{(N+1)}}\beta_h\varphi({\cal A}x\dot-({\cal A}\tilde h\dot+
h)).
$$
Since  $\tilde h\in H_0^{(m)}$ then ${\cal A}\tilde h\in H_0$.
Therefore ${\cal A}\tilde h\dot+h\in H_0$. This means that $f\in
{\rm span}(\varphi({\cal A}x\dot-h))_{h\in H_0}$. It follows
$V_0\subset V_1$. $\square$

 \begin{Th}
  Let $(\varphi(x\dot-h))_{h\in H_0}$ be an orthonormal basis  in
 $V_0$.  Then  $\bigcap\limits_{n\in \mathbb Z}V_n=\{0\}$.
 \end{Th}
 {\bf Proof.} Let $f\in V_{-n}$ for some  $n\in\mathbb N$. Then
 $f({\cal A}^nx)\in V_0$. Since the system   $(\varphi(x\dot-h))_{h\in H_0}$
 is orthonormal we have the equality
 $$
 \frac{1}{p^n}\sum_{h\in H_0}\left|\int\limits_{\mathfrak G}f(x)\varphi({\cal
 A}^{-n}x\dot-h)\,d\mu\right|^2=\sum_{h\in
 H_0}\left|\int\limits_{\mathfrak G}f({\cal
 A}^nx)\varphi(x\dot-h)\,d\mu\right|^2=
 $$
 $$
 =||f({\cal A}^nx)||_2^2=\int\limits_{\mathfrak G}|f({\cal
 A}^nx)|^2\,d\mu=\frac{1}{p^n}||f||_2^2.
 $$
  It is evident that $(p^{\frac{n}{2}}\varphi({\cal
 A}^nx\dot-h))_{h\in H_0}$ is orthonormal basis in $V_n$.
 Therefore
 $$
 ||f||_2^2=p^n\sum_{h\in H_0}\left|\int\limits_\mathfrak Gf(x)\varphi({\cal
 A}^{-n}x\dot-h)\,d\mu\right|
 $$
 for $f\in V_n$.
 Combining  these equations  we obtain
  $$
 ||f||_2^2=\frac{1}{p^n}\sum_{h\in
 H_0}\left|\int\limits_{\mathfrak G}f(x)\varphi({\cal
 A}^{-n}x\dot-h)\,d\mu\right|^2=\frac{1}{p^n}||f||_2^2,
 $$
 for any $n\in \mathbb N$. It follows  $f(x)=0$ a.e. $\square$

 \begin{Th}
  Let $\varphi$ be a solution of the  equation  (\ref{eq3.2}) and  $(\varphi(x\dot-h))_{h\in H_0}$ an orthonormal basis  in
 $V_0$.
  Then $\overline{\bigcup\limits_{n\in\mathbb Z}V_n}=L_2(\mathfrak G)$ if
and only if
 $$
   \bigcup\limits_{n\in\mathbb Z}{\rm supp}\,\hat\varphi(\cdot
   {\cal A}^{-n})=X.
 $$
 \end{Th}
 {\bf Proof.} This theorem is written in \cite{14} for any zero-dimensional group under the condition
 $|\hat{\varphi}|={\bf  1}_{\mathfrak G_0^{\perp}}$.
 But this condition was used to get the inclusion $V_n\subset
 V_{n+1}$ only. By theorems 3.2 the inclusion $V_n\subset
 V_{n+1}$ holds.  Therefore the theorem is true.   $\square$

 The refinement equation (\ref{eq3.2}) may be written in the form
 \begin{equation}                                      \label{eq3.4}
 \hat\varphi(\chi)=m_0(\chi)\hat\varphi(\chi{\cal
  A}^{-1}),
 \end{equation}
  where

 \begin{equation}                                      \label{eq3.5}
 m_0(\chi)=\frac{1}{p}\sum_{h\in
 H_0^{(N+1)}}\beta_h\overline{(\chi{\cal A}^{-1},h)}
 \end{equation}
is a mask of the equation (\ref{eq3.4}).
\begin{Lm}
Let $\varphi\in\mathfrak D_M(\mathfrak G_{-N})$. Then the mask
$m_0(\chi)$ is constant on cosets $\mathfrak G_{-N}^\bot\zeta$.
\end{Lm}
{\bf Proof.} We will prove that $(\chi,{\cal A}^{-1}h)$ are
constant on cosets $\mathfrak G_{-N}^\bot\zeta$. Without loss of
generality, we can assume that $\zeta=r_{-N}^{\alpha_{-N}}\dots
r_{-N+s}^{\alpha_{-N}+s}\notin \mathfrak G_{-N}^\bot$. If
$$
h=a_{-1}g_{-1}\dot+\dots\dot+a_{-N-1}g_{-N-1}\in H_0^{(N+1)}
$$
 then
$$
{\cal A}^{-1}h=a_{-1}g_{0}\dot+\dots\dot+a_{-N-1}g_{-N}\in
\mathfrak G_{-N}.
$$
 If $\chi\in \mathfrak G_{-N}^\bot\zeta$ then
$\chi=\chi_{-N}\zeta$ where  $\chi_{-N}\in \mathfrak G_{-N}^\bot$.
Therefore $(\chi,{\cal A}^{-1}h)=(\chi_{-N}\zeta,{\cal
A}^{-1}h)=(\zeta,{\cal A}^{-1}h)$. This means that $(\chi,{\cal
A}^{-1}h)$ depends on  $\zeta$ only. $\square$

 \begin{Lm}
 The mask  $m_0(\chi)$ is a periodic  function with any period
 $r_1^{\alpha_1}r_2^{\alpha_2}\dots r_s^{\alpha_s}$ $(s\in\mathbb
 N,\; \alpha_j=\overline{0,p-1},\;j=\overline{1,s})$.
 \end{Lm}
 {\bf Proof.}
 Using the equation $(r_k,g_l)=1, (k\neq l)$ we find
$$
(\chi r_1^{\alpha_1}r_2^{\alpha_2}\dots r_s^{\alpha_s},{\cal
A}^{-1}h)=(\chi r_1^{\alpha_1}r_2^{\alpha_2}\dots r_s^{\alpha_s},
a_{-1}g_0\dot+a_{-2}g_{-1}\dot+\dots\dot+a_{-N-1}g_{-N})=
$$
$$
=(\chi,a_{-1}g_0\dot+a_{-2}g_{-1}\dot+\dots\dot+a_{-N-1}g_{-N})=(\chi{\cal
A}^{-1},h).
$$
Therefore  $m_0(\chi r_1^{\alpha_1}\dots
r_s^{\alpha_s})=m_0(\chi)$ and the lemma  is proved. $\square$
\begin{Lm}
The mask $m_0(\chi)$ is defined by its values on cosets $\mathfrak
G_{-N}^\bot r_{-N}^{\alpha_{-N}}\dots r_0^{\alpha_0}$
$(\alpha_j=\overline{0,p-1})$.
\end{Lm}
{\bf Proof.} Let us denote
$$
k=\alpha_0+\alpha_{-1}p+\dots+\alpha_{-N}p^N\in [0,p^{N+1}-1],
$$
$$
l=a_{-1}+a_{-2}p+\dots+a_{-N-1}p^N\in [0,p^{N+1}-1].
$$
Then (\ref{eq3.5}) can be written as the system
\begin{equation} \label{eq3.6}
m_0(\chi_k)=\frac{1}{p}\sum_{l=0}^{p^{N+1}-1}\beta_l\overline{(\chi_k,{\cal
A}^{-1}h_l)},\;k=\overline{0,p^{N+1}-1}
 \end{equation}
in the unknowns $\beta_l$. We consider the characters $\chi_k$ on
the subgroup $\mathfrak G_{-N_0}$. Since ${\cal A}^{-1}h_l$ lie in
$\mathfrak G_{-N}$, it follows that the matrix
$p^{-\frac{N+1}{2}}\overline{(\chi_k,{\cal A}^{-1}h_l)}$ is
unitary, and so the system (\ref{eq3.6}) has a unique solution for
each finite sequence $(m_0(\chi_k))_{k=0}^{p^{N+1}-1}$.
\rightline{$\square$}
 {\bf Remark.} The function $m_0(\chi)$
constructing in  Lemma 3.7 may be not a mask for $\varphi\in
\mathfrak D_M(\mathfrak G_{-N})$. In the section 4 we find
conditions under which the function $m_0(\chi)$ will be a mask.
\begin{Lm}
Let $\hat f_0(\chi)\in{\mathfrak D}_{-N}(\mathfrak G_1^\bot)$.
Then
\begin{equation} \label{eq3.7}
\hat f_0(\chi)=\frac{1}{p}\sum_{h\in
H_0^{(N+1)}}\beta_h\overline{(\chi,{\cal A}^{-1}h)}.
 \end{equation}
\end{Lm}
{\bf Prof.} Since $\int\limits_{\mathfrak
G_0^\bot}(\chi,g)\overline{(\chi,h)}\,d\nu(\chi)=\delta_{h,g}$
for $h,g\in H_0$ it follows that\\
$\int\limits_{\mathfrak G_0^\bot}(\chi{\cal
A}^{-1},g)\overline{(\chi{\cal
A}^{-1},h)}\,d\nu(\chi)=p\delta_{h,g}$.\\
 Therefore we can
consider the set $\left(\frac{{\cal
A}^{-1}h}{\sqrt{p}}\right)_{h\in H_0^{(N+1)}}$ as an orthonormal
system on $\mathfrak G_1^\bot$. We know (lemma 3.5) that
$(\chi,{\cal A}^{-1}h)$ is a constant on cosets $\mathfrak
G_{-N}^\bot\zeta$. It is evident the dimensional of $\mathfrak
D_{-N}(\mathfrak G_1^\bot)$ is equal to $p^{N+1}$. Therefore the
system $\left(\frac{{\cal A}^{-1}h}{\sqrt{p}}\right)_{h\in
H_0^{(N+1)}}$ is an orthonormal basis for $\mathfrak
D_{-N}(\mathfrak G_1^\bot)$ and the equation (\ref{eq3.7}) is
valid. $\square$

%\end{document}
%\begin{Lm}
%Пусть $\varphi$ -- гладкая финитная функция, ${\rm
%supp}\,\hat\varphi(\chi)\subset G_M^\bot$, $\hat\varphi(\chi)$
%постоянна на смежных классах $G_{-N}^\bot\zeta$. Функция $\varphi$
%будет решением некоторого масштабирующего уравнения
%$\hat\varphi(\chi)=m_0(\chi)\hat\varphi(\chi{\cal A}^{-1})$ с
%функцией $m_0(\chi)\in {\mathfrak D}_{-N}(G_m^\bot)$ тогда и
%только тогда, когда ${\rm supp}\,\hat\varphi(\chi)\subset {\rm
%supp}\,\hat\varphi(\chi{\cal A}^{-1})$.
%\end{Lm}
%{\bf Доказательство.} Пусть $\hat\varphi(\chi)$ -- решение
%уравнения $\hat\varphi(\chi)=m_0(\chi)\hat\varphi(\chi{\cal
%A}^{-1})$ и пусть $\hat\varphi(\chi)\ne 0$ на некотором смежном
%классе $G_{-N}^\bot\zeta$. Тогда на этом же смежном классе
%$\hat\varphi(\chi{\cal A}^{-1})\ne 0$. Это и означает, что ${\rm
%supp}\,\hat\varphi(\chi)\subset {\rm supp}\,\hat\varphi(\chi{\cal
%A}^{-1})$.
%
%Обратно, пусть ${\rm supp}\,\hat\varphi(\chi)\subset {\rm
%supp}\,\hat\varphi(\chi{\cal A}^{-1})$. Если на смежном классе
%$G_{-N}^\bot\zeta$ преобразование Фурье $\hat\varphi(\chi)\ne 0$,
%то $\hat\varphi(\chi{\cal A}^{-1})\ne 0$ на этом же смежном
%классе. Поэтому положим
%$m_0(\chi)=\frac{\hat\varphi(\chi)}{\hat\varphi(\chi{\cal
%A}^{-1})}$ на этом смежном классе. Если
%$\hat\varphi(G_{-N}^\bot\zeta)=0$, то получаем
%$m_0(G_{-N}^\bot\zeta)=0$. Тогда выполнено равенство
%$\hat\varphi(\chi)=m_0(\chi)\hat\varphi(\chi{\cal A}^{-1})$ на
%любом смежном классе $G_{-N}^\bot\zeta$. $\square$
%
%
%
 \section{The main results. The statements and proofs}
 In this section we find the necessary and sufficient condition
 under which a step function $\varphi(x)\in {\mathfrak
 D}_M(\mathfrak G_{-N})$ generates an orthogonal MRA on the $p$-adic Vilenkin
 group.
   We  will prove also  that for any $n\in \mathbb N$ there exists a step function $\varphi$
  such that 1) $\varphi$ generate   an orthogonal MRA, 2) ${\rm supp}\,\hat{\varphi}\subset
 \mathfrak G_n^{\perp}$, 3) $\hat\varphi(\mathfrak G_{n}^\bot
  \setminus \mathfrak
G_{n-1}^\bot)\not\equiv 0$.

 First we obtain a test under which the
system of
  shifts $(\varphi(x\dot-h))_{h\in H_0}$ is an orthonormal system.
 \begin{Th}
 Let  $\varphi(x)\in {\mathfrak
 D}_M(\mathfrak G_{-N})$. A shift's system
 $(\varphi(x\dot-h))_{h\in H_0}$ will be orthonormal if and only
if for any
$\alpha_{-N},\alpha_{-N+1},\dots,\alpha_{-1}=\overline{(0,p-1)}$
 \begin{equation} \label{eq4.1}
 \sum_{\alpha_{0},\alpha_1,\dots,\alpha_{M-1}=0}^{p-1}|\hat\varphi(\mathfrak G_{-N}^\bot
 r_{-N}^{\alpha_{-N}}\dots r_0^{\alpha_0}\dots
 r_{M-1}^{\alpha_{M-1}})|^2=1.
 \end{equation}
 \end{Th}
 {\bf Proof.} First we prove that the system
 $(\varphi(x\dot-h))_{h\in H_0}$ will be orthonormal if and only
 if
 \begin{equation} \label{eq4.2}
 \sum_{\alpha_{-N},\dots,\alpha_0,\dots,\alpha_{M-1}}|\hat\varphi(\mathfrak G_{-N}^\bot
 r_{-N}^{\alpha_{-N}}\dots r_{M-1}^{\alpha_{M-1}})|^2=p^N.
 \end{equation}
 and for any vector
 $(a_{-1},a_{-2},\dots,a_{-N}) \ne(0,0,\dots,0),\ (a_j=0,p-1)$
 \begin{equation} \label{eq4.3}
 \sum_{\alpha_{-1},\dots,\alpha_{-N}}\exp\left(\frac{2\pi
 i}{p}(a_{-1}\alpha_{-1}+a_{-2}\alpha_{-2}+\dots+a_{-N}\alpha_{-N})\right)\times
$$
  $$
  \times\sum_{\alpha_0,\alpha_1,\dots,\alpha_{M-1}}|\hat\varphi(\mathfrak G_{-N}^\bot
 r_{-N}^{\alpha_{-N}}\dots r_{M-1}^{\alpha_{M-1}})|^2=0
 \end{equation}

 Let $(\varphi(x\dot-h))_{h\in H_0}$  be an orthonormal system. Using  the Plansherel equality and
 Lemma  2.3 we have
 $$
 \delta_{h_1h_2}=\int\limits_\mathfrak G\varphi(x\dot-h_1)\overline{\varphi(x\dot-h_2)}\,d\mu(x)=
 \int\limits_{\mathfrak G_M^\bot}|\hat\varphi(\chi)|^2(\chi,h_2\dot-
 h_1)\,d\nu(\chi)=
 $$
 $$
 =\sum_{\alpha_{-N},\dots,\alpha_0,\dots,\alpha_{M-1}}\int\limits_{\mathfrak G_{-N}^\bot
 r_{-N}^{\alpha_{-N}}\dots r_0^{\alpha_0}\dots
 r_{M-1}^{\alpha_{M-1}}}|\hat\varphi(\chi)|^2(\chi,h_2\dot-
 h_1)\,d\nu(\chi)=
 $$
 $$
 =\sum_{\alpha_{-N},\dots,\alpha_{M-1}}|\hat\varphi(G_{-N}^\bot
 r_{-N}^{\alpha_{-N}}\dots r_0^{\alpha_0}\dots
 r_{M-1}^{\alpha_{M-1}}|^2\int\limits_{\mathfrak G_{-N}^\bot
 r_{-N}^{\alpha_{-N}}\dots r_0^{\alpha_0}\dots
 r_{M-1}^{\alpha_{M-1}}}(\chi,h_2\dot- h_1)\,d\nu(\chi)=
 $$
 $$
 =p^{-N}{\bf
 1}_{\mathfrak G_{-N}}(h_2\dot-h_1)\times
 $$
 $$
 \times\sum_{\alpha_{-N},\dots,\alpha_{M-1}}|\hat\varphi(\mathfrak G_{-N}^\bot
 r_{-N}^{\alpha_{-N}}\dots r_0^{\alpha_0}\dots
 r_{M-1}^{\alpha_{M-1}}|^2(r_{-N}^{\alpha_{-N}}\dots
 r_0^{\alpha_0}\dots r_{M-1}^{\alpha_{M-1}},h_2\dot-h_1).
 $$
 If  $h_2=h_1$, we obtain the equality (\ref{eq4.3}). If $h_2\ne
 h_1$ then
 \begin{equation} \label{eq4.4}
 h_2\dot-h_1=a_{-1}g_{-1}\dot+\dots\dot+a_{-N}g_{-N}\in \mathfrak G_{-N}
 \end{equation}
 or
 \begin{equation} \label{eq4.5}
 h_2\dot-h_1=a_{-1}g_{-1}\dot+\dots\dot+a_{-N}g_{-N}\dot+
 \dots\dot+a_{-s}g_{-s}\in \mathfrak G\backslash \mathfrak G_{-N}.
  \end{equation}
  If the condition (\ref{eq4.5}) are fulfilled,  then ${\bf 1}_{\mathfrak G_{-N}}(h_2\dot-h_1)=0$.
  If the condition  (\ref{eq4.4})
 are fulfilled, then
 $$
 {\bf 1}_{\mathfrak G_{-N}}(h_2\dot-h_1)=1,
 $$
 $$
  (r_{-N}^{\alpha_{-N}}\dots
 r_0^{\alpha_0}\dots
 r_{M-1}^{\alpha_{M-1}},h_2\dot-h_1)=(r_{-N},g_{-N})^{a_{-N}\alpha_{-N}}\dots
 (r_{-1},g_{-1})^{a_{-1}\alpha_{-1}}.
 $$
 Using the equality $(r_n,g_n)=e^{\frac{2\pi i}{p}}$ we obtain the equality
 (\ref{eq4.3}).The conversely may be proved by analogy.

 Let as show now if for any vector  $(a_{-1},a_{-2},\dots,a_{-N})\ne(0,0,\dots,0)$ the conditions
  (\ref{eq4.2})  (\ref{eq4.3}) are fulfilled, then for any
   $\alpha_{-N},\alpha_{-N+1},\dots,\alpha_{-1}=\overline{0,p-1}$
 \begin{equation} \label{eq4.6}
  \sum_{\alpha_{0},\alpha_1,\dots,\alpha_{M-1}}|\hat\varphi(\mathfrak G_{-N}^\bot
  r_{-N}^{\alpha_{-N}}\dots r_0^{\alpha_0}\dots
  r_{M-1}^{\alpha_{M-1}})|^2=1.
 \end{equation}
 Let us denote
 $$
  n=\sum_{j=1}^Na_{-j}p^{j-1},\;\;k=\sum_{j=1}^N\alpha_{-j}p^{j-1},\;\;C_{n,k}=e^{\frac{2\pi
  i}{p}\left(\sum_{j=1}^N\alpha_{-j}a_{-j}\right)}.
 $$
 and write the equalities  (\ref{eq4.2}) и (\ref{eq4.3}) as the system
 \begin{equation} \label{eq4.7}
 \begin{array}{l}
  C_{0,0}x_0+C_{0,1}x_1+\dots+C_{0,p^N-1}x_{p^N-1}=p^N \\
  C_{1,0}x_0+C_{1,1}x_1+\dots+C_{1,p^N-1}x_{p^N-1}=0 \\
  \dots \dots  \dots\dots\dots\dots\dots\dots\dots \dots \dots \dots \dots \dots\\
  C_{p^N-1,0}x_0+C_{p^N-1,1}x_1+\dots+C_{p^N-1,p^N-1}x_{p^N-1}=0 \\
 \end{array}
 \end{equation}
  with unknowns
 $$
 x_k=\sum_{\alpha_{0},\alpha_1,\dots,\alpha_{M-1}}|\hat\varphi(G_{-N}^\bot
 r_{-N}^{\alpha_{-N}}\dots r_0^{\alpha_0}\dots
 r_{M-1}^{\alpha_{M-1}})|^2.
 $$
 The matrix $(C_{n,k})$ is orthogonal. Indeed, if\\
 $(a_{-1},a_{-2},\dots,a_{-N})\ne(a'_{-1},a'_{-2},\dots,a'_{-N})$,
 i.e., $k\ne n'$ we obtain
 $$
 \sum_{k=0}^{p^N-1}C_{n,k}\overline{C_{n',k}}=\sum_{\alpha_{-1},\dots,\alpha_{-N}}\exp\left(\frac{2\pi
 i}{p}((a_{-1}-a'_{-1})\alpha_{-1}+(a_{-N}-a'_{-N})\alpha_{-N})\right)=0,
 $$
 so at least one of differences  $a_{-l}-a'_{-l}\ne 0$.
 So, the system (\ref{eq4.7}) has unique solution.
 It is evident that  $x_k=1$ is  a solution of this system. This
 means that
 (\ref{eq4.6}) is fulfil, and the necessity is proved. The
 sufficiency is evident. $\square$

 Now we obtain a necessary and sufficient conditions for function  $m_0(\chi)$ to be
 a mask on the class
 $\mathfrak D_{-N}(\mathfrak G_M^\bot)$, i.e. there exists
 $\hat\varphi\in\mathfrak D_{-N}(\mathfrak G_M^\bot)$ for which
 \begin{equation} \label{eq4.8}
 \hat\varphi(\chi)=m_0(\chi)\hat\varphi(\chi{\cal A}^{-1}).
 \end{equation}
 If $m_0(\chi)$ is a mask of (\ref{eq4.8}) then\\
 T1) $m_0(\chi)$
 is constant on cosets  $\mathfrak G_{-N}^\bot\zeta$,\\
 T2) $m_0(\chi)$ is periodic with any period
  $r_1^{\alpha_1}r_2^{\alpha_2}\dots
 r_s^{\alpha_s}$, $\alpha_j=\overline{0,p-1}$, \\
 T3)
 $m_0(\mathfrak G_{-N}^\bot)=1$. \\
 Therefore we will assume that  $m_0$
 satisfies these conditions.
 Let
 $$
 E_k\subset \mathfrak G_k^\bot\setminus \mathfrak
 G_{k-1}^\bot\;\ ,(k=-N+1,-N+2,\dots,0,1,\dots,M,M+1)
 $$
 be a set,  on which $m_0(E_k)=0$.
 Since  $m_0(\chi)$ is constant on cosets
 $\mathfrak G_{-N}^\bot\zeta$, it follows that $E_k$ is a union of such  cosets or $E_k=\emptyset$.
 \begin{Th}
 $m_0(\chi)$ is a mask of some equation on the class $\mathfrak
 D_{-N}(\mathfrak G_M^\bot)$ if and only if
 \begin{equation} \label{eq4.9}
 \bigcup\limits_{k=-N+1}^{M+1}E_k{\cal
 A}^{M+1-k}=\mathfrak G_{M+1}^\bot\setminus \mathfrak G_M^\bot.
 \end{equation}
 \end{Th}
 {\bf Proof.} Since  $m_0(\chi)=1$  on $\mathfrak G_N$
  it follows that $m_0(\chi{\cal
 A}^{-M-N})=1$ for $\chi\in \mathfrak G_M^\bot$. Therefore  $m_0(\chi)$ will
 be a mask if and only if
 \begin{equation} \label{eq4.10}
 m_0(\chi)m_0(\chi{\cal A}^{-1})\dots m_0(\chi{\cal A}^{-M-N})=0
 \end{equation}
 on $\mathfrak G_{M+1}^\bot\setminus \mathfrak G_M^\bot$. Indeed, if
 (\ref{eq4.10}) is true  we set
 $$
 \hat\varphi(\chi)=\prod\limits_{k=0}^\infty m_0(\chi{\cal
 A}^{-k})\in \mathfrak D_{-N}(\mathfrak G_M^\bot).
 $$
 Then $\hat\varphi(\chi)=m_0(\chi)\hat\varphi(\chi{\cal
 A}^{-1}) $  and
 $$
 m_0(\chi)=\sum_{h\in H_0^{(N+1)}}\beta_h\overline{(\chi{\cal
 A}^{-1},h)}
 $$
 for some  $\beta_h$. Therefore $m_0(\chi)$ is a mask. Inversely
  let $m_0(\chi)$ be a mask, i.e.
 $\hat\varphi(\chi)=m_0(\chi)\hat\varphi(\chi{\cal A}^{-1})\in \mathfrak D_{-N}(\mathfrak G_M^\bot)$.
  From it we find
 $$
 \hat\varphi(\chi)=m_0(\chi)m_0(\chi{\cal A}^{-1})\dots
 m_0(\chi{\cal A}^{-M-N})\hat\varphi(\chi{\cal A}^{-M-N-1}),
 $$
 and  $\hat\varphi(\chi{\cal A}^{-M-N-1})=1$ on $\mathfrak G_{M+1}^\bot$.
 Since  $\hat\varphi(\chi)=0$ on $\mathfrak G_{M+1}^\bot\setminus \mathfrak G_M^\bot$,
 it follows
 $$
 m_0(\chi)m_0(\chi{\cal A}^{-1})\dots m_0(\chi{\cal A}^{-M-N})=0
 $$
 on  $\mathfrak G_{M+1}^\bot\setminus \mathfrak G_M^\bot$.

  To
 conclude the proof, it remains to note that  for any
 $-N+1\le k\le M+1$ the inclusion  $E_k{\cal
 A}^{-k+M+1}\subset \mathfrak G_{M+1}^\bot\setminus \mathfrak G_M^\bot$ is true.
 Therefore the equation (\ref{eq4.9}) is fulfil if and only if the equation (\ref{eq4.10}) is true. $\square$
\begin{Lm}
Let $\hat\varphi\in \mathfrak D_{-N}(\mathfrak G_M^\bot)$ be a
solution of the refinement equation
$$
\hat\varphi(\chi)=m_0(\chi)\hat\varphi(\chi{\cal A}^{-1}).
$$
Then for any
$\alpha_{-N},\alpha_{-N+1},\dots,\alpha_{-1}=\overline{0,p-1}$
 \begin{equation} \label{eq4.11}
\sum_{\alpha_0=0}^{p-1}|m_0(\mathfrak G_{-N}^\bot
r_{-N}^{\alpha_{-N}}r_{-N+1}^{\alpha_{-N+1}}\dots
r_{-1}^{\alpha_{-1}}r_{0}^{\alpha_{0}})|^2=1.
 \end{equation}
\end{Lm}
{\bf Proof.} Since $\hat\varphi\in\mathfrak D_{-N}(\mathfrak
G_M^\bot)$, it follows that $\hat\varphi(\mathfrak
G_{M+1}^\bot\setminus \mathfrak G_M^\bot)=0$. Using  theorem 4.1
we have
$$
1=\sum_{\alpha_0,\alpha_1,\dots,
\alpha_{M-1}=0}|\hat\varphi(\mathfrak G_{-N}^\bot
r_{-N}^{\alpha_{-N}}\dots r_0^{\alpha_0}\dots
r_{M-1}^{\alpha_{M-1}})|^2=
$$
$$
=\sum_{\alpha_0,\dots,
\alpha_{M-1},\alpha_{M}=0}|\hat\varphi(\mathfrak G_{-N}^\bot
r_{-N}^{\alpha_{-N}}\dots r_0^{\alpha_0} \dots
r_{M-1}^{\alpha_{M-1}}r_{M}^{\alpha_{M}})|^2=\sum_{\alpha_0=0}^{p-1}|m_0(\mathfrak
G_{-N}^\bot r_{-N}^{\alpha_{-N}}\dots  r_{0}^{\alpha_{0}})|^2
$$
$$
\cdot\sum_{\alpha_1,\dots,
\alpha_{M-1},\alpha_M=0}|\hat\varphi(\mathfrak G_{-N}^\bot
r_{-N}^{\alpha_{-N}+1}\dots r_{-1}^{\alpha_{0}}
r_0^{\alpha_1}\dots r_{M-2}^{\alpha_{M-1}}
r_{M-1}^{\alpha_{M}})|^2=
$$
$$
=\sum_{\alpha_0=0}^{p-1}|\mathfrak G_{-N}^\bot
r_{-N}^{\alpha_{-N}}\dots  r_0^{\alpha_0}|^2.\;\;\square
$$
 {\bf Corollary.} {\it If $N=1$ and $m_0(\mathfrak G_{-1}^\bot
r_{-1}^{\alpha_{-1}}r_{0}^{\alpha_{0}})=\lambda_{\alpha_{-1}+\alpha_0p}$
then we can write the equations (\ref{eq4.11}) in the form
 \begin{equation} \label{eq4.12}
\sum_{\alpha_0=0}^{p-1}|\lambda_{\alpha_{-1}+\alpha_0p}|^2=1.
 \end{equation}
}
\begin{Th}
Suppose the function $m_0(\chi)$ satisfies the conditions
T1,T2,T3, (\ref{eq4.10}), and the function
$$
\hat\varphi(\chi)=\prod\limits_{n=0}^\infty m_0(\chi{\cal A}^{-n})
$$
satisfies the condition  (\ref{eq4.1}). Then $\varphi\in \mathfrak
D_M(\mathfrak G_{-N})$ generates an orthogonal MRA.
\end{Th}
{\bf Proof.} It is evident that $\hat \varphi\in \mathfrak
D_{-N}(\mathfrak G_{M}^\bot)$,
$\hat\varphi(\chi)=m_0(\chi)\hat\varphi(\chi{\cal A}^{-1})$ and
$(\varphi(x\dot-h))_{h\in H_0}$ is an orthonormal system. From
theorems 3.4, 3.3, 3.2 we find that the function $\varphi$
generates an orthogonal MRA. $\square$
\begin{df}
A mask $m_0(\chi)$ is called $N$-elementary $(N\in\mathbb N)$ if
it is constant on cosets $\mathfrak G_{-N}^\bot\chi$ and its
modulus  $|m_0(\chi)|$ take two values:0 and 1 only. The refinable
function $\varphi$ with Fourier transform
$$
\hat\varphi(\chi)=\prod\limits_{j=0}^\infty m_0(\chi{\cal A}^{-1})
$$
is called $N$-elementary too.
\end{df}
\begin{Th}
Let $m_0(\chi)$ be an $1$-elementary mask such that
$$
\sum\limits_{\alpha_0=0}^{p-1}|m_0(\mathfrak G_{-1}^\bot
r_{-1}^{\alpha_{-1}}r_{0}^{\alpha_{0}})|^2=1
$$
for any $\alpha_{-1}=\overline{0,p-1}$.
 Let us denote
 $$
 E_0^{(0)}=\{\alpha=\overline{0,p-1}:\;m_0(\mathfrak
G_{-1}^\bot r_{-1}^{\alpha})=0\}$$ and $l=\sharp E_0^{(0)} $,
 $0\le
l\le p-2$. If $\hat\varphi(\chi)=\prod\limits_{j=0}^\infty
m_0(\chi{\cal A}^{-j})$, then $\hat\varphi(\mathfrak
G_{l+1}^\bot\setminus\mathfrak G_{l}^\bot)=0$.
\end{Th}
{\bf Proof.} Since
$$
\mathfrak G_{l+1}^\bot\setminus\mathfrak
G_{l}^\bot=\bigsqcup\limits_{\alpha_l=1}^{p-1}\bigsqcup\limits_{\alpha_{l-1},\dots,\alpha_{-1}=0}^{p-1}(\mathfrak
G_{-1}^\bot r_{-1}^{\alpha_{-1}}r_{0}^{\alpha_{0}}\dots
r_{l-1}^{\alpha_{l-1}}r_{l}^{\alpha_{l}})
$$
we need prove that
$$
\hat\varphi(\mathfrak G_{-1}^\bot
r_{-1}^{\alpha_{-1}}r_{0}^{\alpha_{0}}\dots
r_{l-1}^{\alpha_{l-1}}r_{l}^{\alpha_{l}})=0
$$
for $\alpha_l=\overline{1,p-1}$;
$\alpha_{-1},\dots,\alpha_{l-1}=\overline{0,p-1}$. Using a
periodicity of $\varphi$ we can write
$$
\hat\varphi(\mathfrak G_{-1}^\bot
r_{-1}^{\alpha_{-1}}r_{0}^{\alpha_{0}}\dots
r_{l-1}^{\alpha_{l-1}}r_{l}^{\alpha_{l}})=
$$
$$
m_0(\mathfrak G_{-1}^\bot
r_{-1}^{\alpha_{-1}}r_{0}^{\alpha_{0}}\dots
r_{l-1}^{\alpha_{l-1}}r_{l}^{\alpha_{l}})\hat\varphi(\mathfrak
G_{-2}^\bot r_{-2}^{\alpha_{-1}}r_{-1}^{\alpha_{0}}\dots
r_{l-2}^{\alpha_{l-1}}r_{l-1}^{\alpha_{l}})=
$$
$$
m_0(\mathfrak G_{-1}^\bot
r_{-1}^{\alpha_{-1}}r_{0}^{\alpha_{0}})\hat\varphi(\mathfrak
G_{-1}^\bot r_{-1}^{\alpha_{0}}r_{0}^{\alpha_{1}}\dots
r_{l-2}^{\alpha_{l-1}}r_{l-1}^{\alpha_{l}})=\dots=
$$
$$
m_0(\mathfrak G_{-1}^\bot
r_{-1}^{\alpha_{-1}}r_{0}^{\alpha_{0}})m_0(\mathfrak G_{-1}^\bot
r_{-1}^{\alpha_0\alpha_{1}})\dots m_0(\mathfrak G_{-1}^\bot
r_{-1}^{\alpha_{l-1}}r_{0}^{\alpha_{l}})m_0(\mathfrak G_{-1}^\bot
r_{-1}^{\alpha_{l}}).
$$
Let us denote $m_0(\mathfrak G_{-1}^\bot
r_{-1}^{k}r_{0}^{j})=\lambda_{k+jp}$ and write $\hat\varphi$ in
the form
$$
\hat\varphi(\mathfrak G_{-1}^\bot
r_{-1}^{\alpha_{-1}}r_{0}^{\alpha_{0}}\dots
r_{l-1}^{\alpha_{l-1}}r_{l}^{\alpha_{l}})=\lambda_{\alpha_{-1}+\alpha_0p}\cdot\lambda_{\alpha_{0}+\alpha_1p}\dots\lambda_{\alpha_{l-1}+\alpha_lp}
\cdot\lambda_{\alpha_l}.
$$
 We will consider numbers $\lambda_{k+jp}$
as elements of the matrix $\Lambda=(\lambda_{j,k})$, where $j$ is
a number of a line, $k$ is a number of a column. Let us consider
the product
$$
\Pi=\lambda_{\alpha_{-1}+\alpha_0p}\cdot\lambda_{\alpha_{0}+\alpha_1p}\dots\lambda_{\alpha_{l-2}+\alpha_{l-1}p}
\cdot\lambda_{\alpha_{l-1}+\alpha_{l}p}
\cdot\lambda_{\alpha_{l}}\;\;(\alpha_{l}\ne 0).
$$
We need prove that $\Pi=0$ for $\alpha_j=\overline{0,p-1}$,
$j=\overline{-1,l-1}$ and for $\alpha_{l}=\overline{1,p-1}$.

If $\alpha_{l}\in E_0^{(0)}$, then $\lambda_{\alpha_{l}}=0$ and
$\Pi=0$. Let $\alpha_{l}\in E_0^{(1)}$ and $\alpha_{l}\ne 0$.

If $\lambda_{\alpha_{l-1}+\alpha_lp}=0$, then $\Pi=0$ and theorem
is proved. Therefore we assume
$|\lambda_{\alpha_{l-1}+\alpha_lp}|=1$. In this case
$\alpha_{l-1}\in E_0^{(0)}$ and $\alpha_{l-1}=0$. Let us consider
$\lambda_{\alpha_{l-2}+\alpha_{l-1}p}$. If
$\lambda_{\alpha_{l-2}+\alpha_{l-1}p}=0$ then $\Pi=0$. Therefore
we assume $|\lambda_{\alpha_{l-2}+\alpha_{l-1}p}|=1$. In this case
$\alpha_{l-1}\in E_0^{(0)}$ and $\alpha_{l-1}\ne \alpha_{l-1}$.
Let us consider $\lambda_{\alpha_{l-3}+\alpha_{l-2}p}$. If
$\lambda_{\alpha_{l-3}+\alpha_{l-2}p}=0$ then $\Pi=0$ and the
theorem is proved. Therefore we assume
$|\lambda_{\alpha_{l-3}+\alpha_{l-2}p}|=1$. In this case
$\alpha_{l-3}\in E_0^{(0)}$ and $\alpha_{l-3}\notin \{
\alpha_{l-1},\alpha_{l-2}\}$.

In the general case, if
$$
|\lambda_{\alpha_{l-s}+\alpha_{l-s+1}p}|\cdot|\lambda_{\alpha_{l-s+1}+\alpha_{l-s+2}p}|\dots|\lambda_{\alpha_{l-1}+\alpha_{l}p}|\cdot|\lambda_{\alpha_l}|=1
$$
and
$$
\alpha_{l-s}\notin
\{\alpha_{l-s+1},\alpha_{l-s+2},\dots,\alpha_{l-1}\},\;\;\alpha_{l-s}\in
E_0^{(0)}
$$
then we consider $\lambda_{\alpha_{l-s-1}+\alpha_{l-s}p}$. If
$|\lambda_{\alpha_{l-s-1}+\alpha_{l-s}p}|=0$ then $\Pi=0$ and the
theorem is proved. If $|\lambda_{\alpha_{l-s-1}+\alpha_{l-s}p}|=1$
then
$$
\alpha_{l-s-1}\notin
\{\alpha_{l-s},\alpha_{l-s+1},\dots,\alpha_{l-1}\},\;\;\alpha_{l-s-1}\in
E_0^{(0)}.
$$
We have two possible cases.\\
1) For some $s\le l$
$$
\lambda_{\alpha_{l-s}+\alpha_{l-s+1}p}\cdot\lambda_{\alpha_{l-s+1}+\alpha_{l-s+2}p}\dots\lambda_{\alpha_{l-1}+\alpha_{l}p}\cdot\lambda_{\alpha_l}=0.
$$
In this case $\Pi=0$, and the theorem is proved.\\
2) For $s=l$
$$
|\lambda_{\alpha_{0}+\alpha_{1}p}|\cdot|\lambda_{\alpha_{1}+\alpha_{2}p}|\dots|\lambda_{\alpha_{l-1}+\alpha_{l}p}|\cdot|\lambda_{\alpha_l}|=1.
$$
In this case $\lambda_{\alpha_{-1}+\alpha_0p}=0$ for
$\alpha_{-l}=\overline{0,p-1}$, then $\Pi=0$ and the theorem is
proved. $\square$

{\bf Remark.} If $l=p-1$, then $m_0(\mathfrak G_0^\bot\setminus
\mathfrak G_{-1}^\bot)\equiv 0$. It follow $\hat\varphi(\mathfrak
G_0^\bot\setminus \mathfrak G_{-1}^\bot)$ and consequently ${\rm
supp}\,\hat\varphi(\chi)=\mathfrak G_{-1}^\bot$. In this case the
system of shift $(\varphi(x\dot-h))_{h\in H_0}$ is not orthonormal
system.

 If $l=0$, then $|m_0(\mathfrak G_0^\bot)|\equiv 1$ and the system
 of shifts $(\varphi(x\dot-h))_{h\in H_0}$ will be orthonormal if
 and only if $\hat\varphi(\mathfrak G_1^\bot\setminus \mathfrak
 G_{0}^\bot)\equiv 0$. In this case $\varphi$ generate an orthogonal MRA on any zero-dimensional group [14].  \\
 {\bf Corollary.} {\it Let $\varphi\in\mathfrak D_M(\mathfrak
 G_{-N})$ be an $1$-elementary refinable function and $\varphi$
 generate an orthogonal MRA on $p$-adic Vilenkin group $\mathfrak
 G$ with $p\ge 3$. Then ${\rm supp}\,\hat\varphi(\chi)\subset
 \mathfrak G_{p-2}^\bot$. }

The next theorem shows the sharpness of this result.

\begin{Th}
Let $\mathfrak G$ -- be a $p$-adic Vilenkin group, $p\ge 3$. Then
for any $1\le l\le p-2$ there exists an $1$-elementary refinable
function $\varphi\in {\mathfrak D_l(\mathfrak G_{-1})}$ that
generate an orthogonal MRA on group $\mathfrak G$.
\end{Th}
{\bf Proof.} We will find the Fourier transform $\hat\varphi$ as
product
$$
\hat\varphi(\chi)=\prod\limits_{j=0}^\infty m_0(\chi{\cal
A}^{-j}),
$$
where the $1$-elementary mask $m_0(\chi)$ is constant on cosets
$\mathfrak G_{-1}^\bot r_{-1}^{\alpha_{-1}}r_{0}^{\alpha_{0}}\dots
r_{s}^{\alpha_{s}}$ $(s\in \mathbb N\bigsqcup \{0\})$.  We will
construct the mask $m_0(\chi)$ on the subgroup $\mathfrak
G_1^\bot$ only, since $m_0(\mathfrak G_{-1}^\bot
r_{-1}^{\alpha_{-1}}r_{0}^{\alpha_{0}}\dots
r_{s}^{\alpha_{s}})=m_0(\mathfrak G_{-1}^\bot
r_{-1}^{\alpha_{-1}}r_{0}^{\alpha_{0}})$. We will assume also that
for any $\alpha_{-1}=\overline{0,p-1}$
\begin{equation} \label{eq4.13}
\sum_{\alpha_0=0}^{p-1}|m_0(\mathfrak G_{-1}^\bot
r_{-1}^{\alpha_{-1}}r_{0}^{\alpha_{0}})|^2=1,
 \end{equation}
since this condition is necessary for mask $m_0(\chi)$.

Choose an arbitrary set $E_l^{(0)}\subset\{1,2,\dots,p-1\}$ of
cardinality $\sharp E_l^{(0)}=l$. Let us denote
$E_l^{(1)}=\{1,2,\dots,p-1\}\setminus E_l^{(0)}$ and
$m_0(\mathfrak G_{-1}^\bot
r_{-1}^{\alpha_{-1}}r_{0}^{\alpha_{0}})=\lambda_{\alpha_{-1}+\alpha_0p}$.
First we set
$$
\lambda_0=1,\;\; |\lambda_\alpha|=\left\{\begin{array}{ll}
  0, & \alpha\in E_l^{(0)}, \\
  1, & \alpha\in E_l^{(1)} \\
\end{array}\right..
$$
 Now we will define
$\lambda_{\alpha_{-1}+\alpha_0p}$ for $\alpha_0\ge 1$. It follow
from (\ref{eq4.13}) that $\lambda_{\alpha_{-1}+\alpha_0p}=0$ for
$\alpha_{-1}\in E_l^{(1)}$, $\alpha_0\ge 1$. Choose an arbitrary
$\alpha_{l-1}^{(0)}\in E_l^{(1)}$ and fix it. Now we choose
$\alpha_{l-2}^{(0)}\in E_l^{(0)}$ and  set
$$
|\lambda_{\alpha_{l-2}^{(0)}+\alpha_{l-1}^{(0)}p}|=1,\;\;|\lambda_{\alpha_{l-2}^{(0)}+\alpha
p}|=0\;\mbox{if}\;\alpha\ne\alpha_{l-1}^{(0)}.
$$
If numbers $\alpha_{l-2}^{(0)},\dots,\alpha_{s}^{(0)}\in
E_l^{(0)}$ $(s=l-1,l-2,\dots,0)$ have been choosen we choose
$\alpha_{s-1}^{(0)}\in
E_l^{(0)}\setminus\{\alpha_{l-2}^{(0)},\dots,\alpha_{s}^{(0)}\}$
and set
$$
|\lambda_{\alpha_{s-1}^{(0)}+\alpha_{s}^{(0)}p}|=1,\;\;|\lambda_{\alpha_{s-1}^{(0)}+\alpha
p}|=0\;\mbox{if}\;\alpha\ne\alpha_{s}^{(0)}.
$$
So the mask $m_0(\chi)$ have been defined on the subgroup
$\mathfrak G_1^\bot$ and consequently on the group $\mathfrak G$.

It is evident that
$$
\lambda_{\alpha_{-1}^{(0)}+\alpha_{0}^{(0)}p}\cdot\lambda_{\alpha_{0}^{(0)}+\alpha_{1}^{(0)}p}\dots
\lambda_{\alpha_{l-2}^{(0)}+\alpha_{l-1}^{(0)}p}\cdot\lambda_{\alpha_{l-1}^{(0)}}\ne
0.
$$
Let us show that for any vector
$(\alpha_{-1},\alpha_{0},\dots,\alpha_{l-1})\ne(\alpha_{-1}^{(0)},\alpha_{0}^{(0)},\dots,\alpha_{l-1}^{(0)})$
\begin{equation} \label{eq4.14}
\lambda_{\alpha_{-1}+\alpha_{0}p}\cdot\lambda_{\alpha_{0}+\alpha_{1}p}\dots
\lambda_{\alpha_{l-2}+\alpha_{l-1}p}\cdot\lambda_{\alpha_{l-1}}=0.
 \end{equation}
Indeed, if $\alpha_{l-1}\in E_l^{(0)}$ then
$\lambda_{\alpha_{l-1}}=0$. If $\alpha_{l-1}\in E_l^{(1)}$ and
$\alpha_{l-1}\ne \alpha_{l-1}^{(0)}$   then
$\lambda_{\alpha_{l-2}+\alpha_{l-1}p}=0$. If $\alpha_{l-1}\in
E_l^{(1)}$ and $\alpha_{l-1}=\alpha_{l-1}^{(0)}$ then we denote
$$
s=\min\{j:\;\alpha_j=\alpha_j^{(0)}\}.
$$
For this $s$ we have $\lambda_{\alpha_{s-1}+\alpha_s^{(0)}p}=0$
and the equality (\ref{eq4.14}) is proved. It should be noted that
$\lambda_{\alpha+\alpha_{-1}^{(0)}p}=0$ for
$\alpha=\overline{0,p-1}$. Therefore
\begin{equation} \label{eq4.15}
\lambda_{\alpha+\alpha_{-1}^{(0)}p}\cdot\lambda_{\alpha_{-1}^{(0)}+\alpha_{0}^{(0)}p}\dots
\lambda_{\alpha_{l-2}^{(0)}+\alpha_{l-1}^{(0)}p}\cdot\lambda_{\alpha_{l-1}^{(0)}}=0.
 \end{equation}
Let us show that $\hat\varphi(\mathfrak G_l^\bot\setminus\mathfrak
G_{l-1}^\bot)\not\equiv 0$ and $\hat\varphi(\mathfrak
G_{l+1}^\bot\setminus\mathfrak G_{l}^\bot)\equiv 0$. Since
$m_0(\chi)$ is periodic with any period
$r_1^{\alpha_1}r_2^{\alpha_2}\dots r_s^{\alpha_s}$, it follow that
$$
\hat\varphi(\mathfrak G_{-1}^\bot
r_{-1}^{\alpha_{-1}}r_0^{\alpha_0}\dots r_{l-1}^{\alpha_{l-1}})=
$$
$$
m_0(\mathfrak G_{-1}^\bot r_{-1}^{\alpha_{-1}}r_0^{\alpha_0}\dots
r_{l-1}^{\alpha_{l-1}})m_0(\mathfrak G_{-1}^\bot
r_{-1}^{\alpha_{0}}r_0^{\alpha_1}\dots
r_{l-1}^{\alpha_{l-1}})\dots m_0(\mathfrak G_{-1}^\bot
r_{-1}^{\alpha_{l-2}}r_0^{\alpha_{l-1}})m_0(\mathfrak G_{-1}^\bot
r_{-1}^{\alpha_{l-1}})=
$$
$$
=m_0(\mathfrak G_{-1}^\bot
r_{-1}^{\alpha_{-1}}r_0^{\alpha_0})m_0(\mathfrak G_{-1}^\bot
r_{-1}^{\alpha_{0}}r_0^{\alpha_1})\dots m_0(\mathfrak G_{-1}^\bot
r_{-1}^{\alpha_{l-2}}r_0^{\alpha_{l-1}})m_0(\mathfrak G_{-1}^\bot
r_{-1}^{\alpha_{l-1}})=
$$
$$
=\lambda_{\alpha_{-1}+\alpha_0p}\cdot\lambda_{\alpha_{0}+\alpha_1p}\cdot\dots\cdot\lambda_{\alpha_{l-2}+\alpha_{l-1}p}\cdot\lambda_{\alpha_{l-1}}\ne
0
$$
for $(\alpha_{-1},\alpha_{0},\dots,\alpha_{l-2},\alpha_{l-1})=
(\alpha_{-1}^{(0)},\alpha_{0}^{(0)},\dots,\alpha_{l-2}^{(0)},\alpha_{l-1}^{(0)})$.
This means that $\hat\varphi(\mathfrak G_l^\bot\setminus\mathfrak
G_{l-1}^\bot)\not\equiv 0$.

By analogy
$$
\hat\varphi(\mathfrak G_{-1}^\bot
r_{-1}^{\alpha_{-1}}r_0^{\alpha_0}\dots
r_{l-1}^{\alpha_{l-1}}r_l^{\alpha_l})=\lambda_{\alpha_{-1}+\alpha_0p}\cdot\lambda_{\alpha_{0}+\alpha_1p}\cdot\dots\cdot
\lambda_{\alpha_{l-1}+\alpha_{l}p}\cdot\lambda_{\alpha_{l}}.
$$
If $\alpha_l\in E_l^{(0)}$ then $\lambda_{\alpha_l}=0$. If
$\alpha_l\in E_l^{(1)}$ and $\alpha_l\ne \alpha_{l-1}^{(0)}$ then
$\lambda_{\alpha_{l-1}+\alpha_{l}p}=0$ for any
$\alpha_{l-1}=\overline{0,p-1}$. If $\alpha_l\in E_l^{(1)}$ and
$\alpha_l= \alpha_{l-1}^{(0)}$ we define the number
$$
s=\min\{j:\;\alpha_j=\alpha_{j-1}^{(0)}\}.
$$
Then
$$
\hat\varphi(\mathfrak G_{-1}^\bot
r_{-1}^{\alpha_{-1}}r_0^{\alpha_0}\dots
 r_l^{\alpha_l})=\lambda_{\alpha_{-1}+\alpha_0p}\dots\lambda_{\alpha_{s-1}+\alpha_{s-1}^{(0)}p}
\lambda_{\alpha_{s-1}^{(0)}+\alpha_{s-2}^{(0)}p}\dots
\lambda_{\alpha_{l-2}^{(0)}+\alpha_{l-1}^{(0)}p}\cdot\lambda_{\alpha_{l-1}^{(0)}}=0
$$
since $\lambda_{\alpha_{s-1}+\alpha_{s-1}^{(0)}p}=0$ for any
$\alpha_{s-1}=\overline{0,p-1}$. This means that
$\hat\varphi(\mathfrak G_{l+1}^\bot\setminus\mathfrak
G_{l}^\bot)\equiv 0$. Consequently $\hat\varphi\in{\mathfrak
D}_{-1}(\mathfrak G_l^\bot)$.

Let us show that $(\varphi(x\dot-h))_{h\in H_0}$ is an orthonormal
system. We need show that the sum
$$
S{(\alpha_{-1})}=\sum_{\alpha_0,\alpha_1,\dots,\alpha_{l-1}=0}^{p-1}|\hat\varphi(\mathfrak
G_{-1}^\bot r_{-1}^{\alpha_{-1}}r_0^{\alpha_0}\dots
r_{l-1}^{\alpha_{l-1}})|^2=
$$
$$
=\sum_{\alpha_0,\alpha_1,\dots,\alpha_{l-1}=0}^{p-1}|\lambda_{\alpha_{-1}+\alpha_0p}|^2|\lambda_{\alpha_{0}+\alpha_1p}|^2
\dots
|\lambda_{\alpha_{l-2}+\alpha_{l-1}p}|^2|\lambda_{\alpha_{l-1}}|^2=1
$$
for any $\alpha_{-1}=\overline{0,p-1}$.

Let us consider next possible cases.\\
1) If $\alpha_{-1}=0$ then $\lambda_{\alpha_{-1}+\alpha_0p}\ne 0$
iff $\alpha_0=0$, $\lambda_{\alpha_{0}+\alpha_1p}\ne 0$ iff
$\alpha_1=0$ and so on.

Consequently $S(\alpha_{-1})\ne 0$ iff
$\alpha_{-1}=\alpha_0=\dots=\alpha_{l-1}=0$. It means that
$S(\alpha_{-1})=1$.\\
2) If $\alpha_{-1}\ne 0$ and $\alpha_{-1}\in E_l^{(1)}$ then
$\lambda_{\alpha_{-1}+\alpha_0p}\ne 0$ iff $\alpha_0=0$ and by
analog $S(\alpha_{-1})=1$.\\
3) If $\alpha_{-1}\in E_l^{(0)}$ and
$\alpha_{-1}=\alpha_{-1}^{(0)}$ then
$\lambda_{\alpha_{-1}+\alpha_0p}\ne 0$ iff
$\alpha_0=\alpha_0^{(0)}$,\\
$\lambda_{\alpha_{0}^{(0)}+\alpha_1p}\ne 0$ iff
$\alpha_1=\alpha_1^{(0)}$ and so on. Consequently
$S(\alpha_{-1})\ne 0$ iff $\alpha_0=\alpha_0^{(0)}$,
$\alpha_1=\alpha_1^{(0)},\dots$,
$\alpha_{l-1}=\alpha_{l-1}^{(0)}$. It means that
$S(\alpha_{-1})=1$.\\
4) If $\alpha_{-1}\in E_l^{(0)}$ and
$\alpha_{-1}=\alpha_{j}^{(0)}$ $(j\ge 0)$ then
$\lambda_{\alpha_{-1}+\alpha_0p}\ne 0$ iff
$\alpha_0=\alpha_{j+1}^{(0)}$,\\
$\lambda_{\alpha_{0}+\alpha_1p}\ne 0$ iff
$\alpha_1=\alpha_{j+2}^{(0)}$ and so on,
$\alpha_{l-j-2}=\alpha_{l-1}^{(0)}$.  Then
$\alpha_{l-j-1}=\dots=\alpha_{l-1}=0$. This means that
$S(\alpha_{-1})=1$. $\square$

 By theorem 4.4
$\varphi(x)$ generate an orthogonal MRA. $\square$

\end{document}